\documentclass{amsart}
\theoremstyle{plain}
\newtheorem{Thm}{Theorem}

\newtheorem{Coro}[Thm]{Corollary}
\newtheorem{prop}[Thm]{Proposition}
\newtheorem{Lem}[Thm]{Lemma}
\newtheorem{Def}[Thm]{Definition}
\newtheorem{Rk}[Thm]{Remark}
\newtheorem{Ex}[Thm]{Example}

\begin{document}
\large
\title{Curve Shortening Flow in a Riemannian Manifold}

\author{Li MA and Dezhong Chen}

\address{Department of mathematical sciences \\
Tsinghua university \\
Beijing 100084 \\
China}

\email{lma@math.tsinghua.edu.cn} \dedicatory{}
\date{Oct. 15th, 2003}

\begin{abstract}In this paper, we systemally study the long time behavior of
the curve shortening flow
in a closed or non-compact complete locally Riemannian symmetric
manifold. Assume that we have a global flow. Then we can exhibit a
a limit for the global behavior of the flow. In particular, we
show the following results. 1). Let $\mathbf{M}$ be a compact
locally symmetric space. If the curve shortening flow exists for
infinite time, and
$$
\lim_{t\rightarrow\infty}L(\gamma_{t})>0,
$$
then for every $n>0$,
$$
\lim_{t\rightarrow \infty}\sup(|\frac{D^{n}T}{\partial s^{n}}|)=0.
$$
In particular, the limiting curve exists and is a closed geodesic
in $\mathbf{M}$. 2). For $\gamma_{0}$ is a ramp, we have a global
flow and the flow converges to a geodesic in $C^{\infty}$ norm.
\end{abstract}

 \maketitle

\begin{section}{Introduction}
  In this paper, we systemally study the limiting behavior at infinity
  of the curve shortening flow in a locally symmetric Riemannian manifold. Curve
  shortening flows in the plane, surfaces, and the
  3-dimensional Euclidean space were studied respectively by M.Gage and R.Hamilton
  \cite{GH86}, M.Grayson \cite{G89}, and S.Altschuler and M.Grayson
  \cite{AG91}. Very recently the curve shortening flow
  in a closed Riemannian manifold has been used by G.Perelman \cite{P03}
  to study the Ricci-Hamilton flow. All these papers can be
  considered as models for our study. It is clear that studying
  the curve shortening flow in a general Riemannian manifolds is
  not a easy work, but it is still quite interesting.

   By definition, our curve shortening flow in the compact or complete
    Riemannian
   manifold $(\mathbf{M},g)$ is evolving the initial curve $\gamma_{0}$ along
   the flow
   $$
   \frac{\partial \gamma }{\partial t}=\frac{DT}{\partial s},
   $$
where $T$ is the unit tangent vector of $\gamma_{t}$,
$\frac{DT}{\partial s}$ is the covariant derivative of $T$ in the
direction $T$ in the space $(\mathbf{M},g)$, and $s$ is the
arc-length parameter of $\gamma_{t}$. It is shown that there
always exists a short time flow for the curve shortening problem
in any compact Riemannian manifold (see $\S$ 2 in \cite{GH86}). We
remark that the existence of a short time flow for the curve
shortening problem in a non-compact complete Riemannian manifold
can also be obtained as in \cite{Shi89}.

    In this work, we mainly prove the following

\begin{Thm}[Main Theorem]. Let $\mathbf{M}$ be a compact locally
symmetric space. If the curve shortening flow exists for infinite
time, and
$$
\lim_{t\rightarrow \infty}L(\gamma_{t})>0,
$$
then for every $n>0$,
$$\lim_{t\rightarrow \infty}\sup (|\frac{D^{n}T}{\partial
s^{n}}|)=0.$$ In particular, the limiting curve exists and is a
closed geodesic in $\mathbf{M}$.

2). For $\gamma_{0}$ is a ramp, we have a global flow and the flow
converges to a geodesic in $C^{\infty}$ norm.
\end{Thm}

Roughly speaking, a {\it ramp} is a curve with non-trivial height
with respect to a conformal Killing or just Killing vector field.
One may see section 6 below for concrete definition. This result
is a special case of our Theorem 10 and Theorem 19 below.
   We also obtain many other results for the curve shortening
   problem
   in varied cases, in particular, when $\mathbf{M}$ is the space forms.
   . The results are exhibited in the sequent
   sections. We remark that in many cases, we do not assume
   $\mathbf{M}$ being compact. We only assume that $(\mathbf{M},g)$
is a complete Riemannian manifold. We give a remark on our
assumption of the global flow. In fact, one can easily see from
our estimate in section 3 that in a nontrivial path homotopy
class, we can always have a global flow. This is also true for
ramps. We will study the finite time blow up of the flow in a
separate paper.

\par The paper is
organized as follows:\par In $\S$ 2, we include some fundamental
formulae for the flow.\par In $\S$ 3, we deduce some precise
estimates for the flow in the locally symmetric Riemannian space
which also reprove the standard long time existence result in
$R^3$ .\par In $\S$ 4, we investigate the limiting behavior of the
curve shortening flow with some necessary assumptions on the
initial curves.\par In $\S$ 5, we compute in details the
evolutions of the curve shortening flows in space forms
$\mathbf{H}^{3}$ and $\mathbf{S}^{3}$ respectively.\par In $\S$ 6,
we study the flow for ramps and use them to find closed
geodesics.\par In $\S$ 7, we analyze the curve shortening flow on
manifolds with time-dependent metrics.
\end{section}

\begin{section}{preliminaries}
On an $n$-dimensional Riemannian manifold $(\mathbf{M},g)$, let
$$
\gamma : \mathbf{S}^{1}\times (a,b)\rightarrow \mathbf{M}
$$
be an evolving immersed curve. Denote by $\gamma_{t}$ the
associated trajectory, i.e.,
$$
\gamma_{t}(\cdot)=:\gamma(\cdot,t).
$$
Then the length of $\gamma_{t}$ is
\begin{equation}
L(\gamma_{t})=\int_{\mathbf{S}^{1}}|\frac{d}{du}\gamma_{t}|du=
\int_{\mathbf{S}^{1}}|\frac{\partial \gamma}{\partial
u}|du=\int_{\mathbf{S}^{1}}vdu,\notag
\end{equation}
where
$$
v =: |\frac{\partial \gamma}{\partial u}|
$$
is the speed. We define the arc-length parameter $s$ by
$$
\frac{\partial}{\partial s}=:\frac{1}{v}\frac{\partial}{\partial
u},
$$
which implies
$$ds = vdu
.$$\par As usual, we denote by $T$ the associated unit tangent
vector, i.e.,
$$
T=:\frac{\partial \gamma}{\partial s}=\frac{1}{v}\frac{\partial
\gamma}{\partial u}.
$$
Then the time derivative of length is
\begin{equation}
\hspace{-32mm}\frac{d}{dt}L(\gamma_{t})=\int_{\mathbf{S}^{1}}\frac{\partial
v }{\partial t}du\notag
\end{equation}
\begin{equation}
\hspace{1mm}=\int_{\mathbf{S}^{1}}<\nabla_{t}\frac{\partial
\gamma}{\partial u},T>du\notag
\end{equation}
\begin{equation}
\hspace{1mm}=\int_{\mathbf{S}^{1}}<\nabla_{u}\frac{\partial
\gamma}{\partial t},T>du\notag
\end{equation}
\begin{equation}
\hspace{32mm}=\int_{\mathbf{S}^{1}}\{\frac{\partial}{\partial
u}<\frac{\partial \gamma}{\partial t},T>-<\frac{\partial
\gamma}{\partial t},\nabla_{u}T>\}du\notag
\end{equation}
\begin{equation}
\hspace{4mm}=-\int_{\mathbf{S}^{1}}<\frac{\partial
\gamma}{\partial t},\nabla_{u}T>du\notag
\end{equation}
\begin{equation}
\hspace{3mm}=-\int_{\mathbf{S}^{1}}<\frac{\partial \gamma
}{\partial t},\frac{DT}{\partial s}>ds.\notag
\end{equation}
If $\gamma$ evolves according to the equation
$$
\frac{\partial \gamma }{\partial t}=\frac{DT}{\partial s},
$$
then we find that
$$
\frac{d}{dt}L(\gamma_{t})=-\int_{\mathbf{S}^{1}}k^{2}ds \leq 0,
$$
where
$$
k^{2} =: |\frac{DT}{\partial s}|^{2}
$$
is the curvature squared. This leads us to give the following
\begin{Def}
For a curve shortening flow, we mean an evolving immersed curve
$\gamma(\cdot,t)$ satisfying the evolution equation
\begin{equation}
\frac{\partial \gamma }{\partial t}=\frac{DT}{\partial s}.\notag
\end{equation}
\end{Def}
\begin{Rk}
Obviously, we can regard $\gamma_{t}(\mathbf{S}^{1})$ as a
1-dimensional sub-manifold of $\mathbf{M}$. With the induced
metric from $\mathbf{M}$, its mean curvature vector field is
\begin{equation}
H=(\nabla_{T}T)^{\bot}.\notag
\end{equation}
Noting that
$$
<T,T> \equiv 1,
$$
we have
$$
<\nabla_{T}T,T> \equiv 0.
$$
This implies
$$
\nabla_{T}T \perp T\gamma_{t}.
$$
So
\begin{equation}
H=\nabla_{T}T.\notag
\end{equation}
This shows that a curve shortening flow is a mean curvature flow.
\end{Rk}
Next we will give some fundamental computations. These formulae
have already appeared in many papers, for example ~\cite{G89}.
\begin{Lem}
The evolution of $v$ is
\begin{equation}
\frac{\partial v}{\partial t}=-k^{2}v.\notag
\end{equation}
\end{Lem}
\begin{proof}
By definition,
$$
v^{2} = <\frac{\partial \gamma}{\partial u},\frac{\partial
\gamma}{\partial u}>.
$$
Differentiating it with respect to $t$, we get
\begin{equation}
\hspace{-12mm}2v\frac{\partial v}{\partial
t}=2<\nabla_{t}\frac{\partial \gamma}{\partial u},\frac{\partial
\gamma}{\partial u}>\notag
\end{equation}
\begin{equation}
\hspace{-2mm}=2<\nabla_{u}\frac{\partial \gamma}{\partial
t},\frac{\partial \gamma}{\partial u}>\notag
\end{equation}
\begin{equation}
\hspace{1mm}=2v^{2}<\nabla_{T}\frac{DT}{\partial s},T>\notag
\end{equation}
\begin{equation}
\hspace{2mm}=-2v^{2}<\frac{DT}{\partial s},\frac{DT}{\partial
s}>\notag
\end{equation}
\begin{equation}
\hspace{-19mm}=-2k^{2}v^{2}\notag
\end{equation}
\end{proof}
\begin{Lem}
Covariant differentiation with respect to $s$ and $t$ are related
by the equation
$$
\nabla_{t}\nabla_{s}=\nabla_{s}\nabla_{t}+k^{2}\nabla_{s}+R(T,\frac{DT}{\partial
s}),
$$
where $R$ is the curvature operator on $\mathbf{M}$.
\end{Lem}
\begin{proof}
We have (see ~\cite{C})
$$
\nabla_{t}\nabla_{u}=\nabla_{u}\nabla_{t}+R(\frac{\partial}{\partial
u},\frac{\partial}{\partial t})\notag
$$
and
\begin{equation}
\nabla_{s}=\frac{1}{v}\nabla_{u}.\notag
\end{equation}
Using Lemma 3, we get
\begin{equation}
\hspace{-15mm}\nabla_{t}\nabla_{s}=\frac{\partial}{\partial
t}(\frac{1}{v})\nabla_{u}+\frac{1}{v}\nabla_{t}\nabla_{u}\notag
\end{equation}
\begin{equation}
\hspace{16mm}=k^{2}\frac{1}{v}\nabla_{u}+\frac{1}{v}\nabla_{u}\nabla_{t}+
\frac{1}{v}R(\frac{\partial}{\partial u},\frac{\partial}{\partial
t})\notag
\end{equation}
\begin{equation}
\hspace{8mm}=\nabla_{s}\nabla_{t}+k^{2}\nabla_{s}+R(T,\frac{DT}{\partial
s}).\notag
\end{equation}
\end{proof}
\begin{Lem}
The covariant differentiation of $T$ with respect to time $t$ is
\begin{equation}
\nabla_{t}T=k^{2}T+\frac{D^{2}T}{\partial s^{2}}.\notag
\end{equation}
\end{Lem}
\begin{proof}
The proof is a straightforward calculation.
\begin{equation}
\hspace{-11mm} \nabla_{t}T=\nabla_{t}(\frac{1}{v}\frac{\partial
\gamma}{\partial u})\notag
\end{equation}
\begin{equation}
\hspace{11mm} =k^{2}\frac{1}{v}\frac{\partial \gamma }{\partial
u}+\frac{1}{v}\nabla_{u}\frac{\partial \gamma }{\partial t}\notag
\end{equation}
\begin{equation}
\hspace{1mm} =k^{2}T+\frac{D^{2}T}{\partial s^{2}}.\notag
\end{equation}
\end{proof}
\end{section}

\begin{section}{Bernstein type estimates }
In this section and the next one, we shall assume that
$\mathbf{M}$ is a locally symmetric space, i.e.,
$$
\nabla R = 0.
$$
For a locally symmetric space, we have
\begin{equation}
\hspace{-57mm}\nabla R(X,Y,Z,W)=(\nabla_{X}R)(Y,Z,W)\notag
\end{equation}
\begin{equation}
\hspace{30mm}=\nabla_{X}(R(Y,Z,W))-R(\nabla_{X}Y,Z,W)-R(Y,\nabla_{X}Z,W)\notag
\end{equation}
\begin{equation}
\hspace{-25mm}-R(Y,Z,\nabla_{X}W)\notag
\end{equation}
\begin{equation}
\hspace{-55mm}=0\notag
\end{equation}
$\Rightarrow$
 \setcounter{equation}{0}
\begin{equation}
\nabla_{X}(R(Y,Z,W))=R(\nabla_{X}Y,Z,W)+R(Y,\nabla_{X}Z,W)
+R(Y,Z,\nabla_{X}W),
\end{equation}
for all $X,Y,Z,W\in T\mathbf{M}$. We shall also assume that
$\mathbf{M}$ satisfies Condition ($\Lambda$), i.e., there exists a
positive constant $\Lambda$, such that
$$
R(\widetilde{X},\widetilde{Y},\widetilde{Z},\widetilde{W})\leq
\Lambda,
$$
for all unit vectors
$\widetilde{X},\widetilde{Y},\widetilde{Z},\widetilde{W}$.\par
With these assumptions, we can give some precise estimates which
will bound the evolution of $|\frac{D^{n}T}{\partial s^{n}}|^{2}$.
This type of estimate also appears in ~\cite{AG91} for the curve
shortening flow in $\mathbf{R}^{3}$. In that case, these estimates
are dilation-invariant, and play an important role in singularity
analysis.\par First, let us compute the time derivative of
$|\frac{D^{n}T}{\partial s^{n}}|^{2}$ as follows:
\begin{equation}
\hspace{-89mm} \frac{\partial}{\partial t}(|\frac{D^{n}T}{\partial
s^{n}}|^{2})\notag
\end{equation}
\begin{equation}
\hspace{-79mm}=2<\frac{D}{\partial t}\frac{D^{n}T}{\partial
s^{n}},\frac{D^{n}T}{\partial s^{n}}>\notag
\end{equation}
\begin{equation}
\hspace{-19mm}=2<\frac{D}{\partial s}\frac{D}{\partial
t}\frac{D^{n-1}T}{\partial s^{n-1}}+k^{2}\frac{D^{n}T}{\partial
s^{n}}+R(T,\frac{DT}{\partial s})\frac{D^{n-1}T}{\partial
s^{n-1}},\frac{D^{n}T}{\partial s^{n}}>\notag
\end{equation}
\begin{equation}
=2<\frac{D}{\partial s}\frac{D}{\partial
t}\frac{D^{n-1}T}{\partial s^{n-1}},\frac{D^{n}T}{\partial
s^{n}}>+2k^{2}|\frac{D^{n}T}{\partial
s^{n}}|^{2}+2R(T,\frac{DT}{\partial s},\frac{D^{n-1}T}{\partial
s^{n-1}},\frac{D^{n}T}{\partial s^{n}})\notag
\end{equation}
\begin{equation}
\hspace{-8mm} =2<\frac{D}{\partial s}(\frac{D}{\partial
s}\frac{D}{\partial t}\frac{D^{n-2}T}{\partial
s^{n-2}}+k^{2}\frac{D^{n-1}T}{\partial
s^{n-1}}+R(T,\frac{DT}{\partial s})\frac{D^{n-2}T}{\partial
s^{n-2}}),\frac{D^{n}T}{\partial s^{n}}>\notag
\end{equation}
\begin{equation}
\hspace{-43mm}+2k^{2}|\frac{D^{n}T}{\partial
s^{n}}|^{2}+2R(T,\frac{DT}{\partial s},\frac{D^{n-1}T}{\partial
s^{n-1}},\frac{D^{n}T}{\partial s^{n}})\notag
\end{equation}
\begin{equation}
\hspace{2mm}=2<\frac{D^{2}}{\partial s^{2}}\frac{D}{\partial
t}\frac{D^{n-2}T}{\partial s^{n-2}},\frac{D^{n}T}{\partial
s^{n}}>+2<\frac{D}{\partial s}(k^{2}\frac{D^{n-1}T}{\partial
s^{n-1}}),\frac{D^{n}T}{\partial
s^{n}}>+2k^{2}|\frac{D^{n}T}{\partial s^{n}}|^{2}\notag
\end{equation}
\begin{equation}
\hspace{-2mm}+2<\frac{D}{\partial s }(R(T,\frac{DT}{\partial
s})\frac{D^{n-2}T}{\partial s^{n-2}}),\frac{D^{n}T}{\partial
s^{n}}>+2R(T,\frac{DT}{\partial s},\frac{D^{n-1}T}{\partial
s^{n-1}},\frac{D^{n}T}{\partial s^{n}})\notag
\end{equation}
\begin{equation}
\hspace{-14mm}=\cdots =2<\frac{D^{n}}{\partial
s^{n}}\frac{D}{\partial t}T,\frac{D^{n}T}{\partial
s^{n}}>+2\sum^{n-1}_{i=0}<\frac{D^{i}}{\partial
s^{i}}(k^{2}\frac{D^{n-i}T}{\partial
s^{n-i}}),\frac{D^{n}T}{\partial s^{n}}>\notag
\end{equation}
\begin{equation}
\hspace{-21mm}+2\sum^{n-1}_{i=0}<\frac{D^{i}}{\partial
s^{i}}(R(T,\frac{DT}{\partial s})\frac{D^{n-1-i}T}{\partial
s^{n-1-i}}),\frac{D^{n}T}{\partial s^{n}}>.\notag
\end{equation}
Using Lemma 5, we get\par
\begin{equation}
\hspace{-39mm}2<\frac{D^{n}}{\partial s^{n}}\frac{D}{\partial
t}T,\frac{D^{n}T}{\partial s^{n}}>\notag
\end{equation}
\begin{equation}
=2<\frac{D^{n+2}T}{\partial s^{n+2}},\frac{D^{n}T}{\partial
s^{n}}>+2<\frac{D^{n}T}{\partial
s^{n}}(k^{2}T),\frac{D^{n}T}{\partial s^{n}}>\notag
\end{equation}
\begin{equation}
\hspace{14mm}=\frac{\partial^{2}}{\partial
s^{2}}(|\frac{D^{n}T}{\partial
s^{n}}|^{2})-2|\frac{D^{n+1}T}{\partial
s^{n+1}}|^{2}+2<\frac{D^{n}T}{\partial
s^{n}}(k^{2}T),\frac{D^{n}T}{\partial s^{n}}>.\notag
\end{equation}
So
\begin{equation}
\frac{\partial}{\partial t}(|\frac{D^{n}T}{\partial
s^{n}}|^{2})=\frac{\partial^{2}}{\partial
s^{2}}(|\frac{D^{n}T}{\partial
s^{n}}|^{2})-2|\frac{D^{n+1}T}{\partial
s^{n+1}}|^{2}+2\sum^{n}_{i=0}<\frac{D^{i}}{\partial
s^{i}}(k^{2}\frac{D^{n-i}T}{\partial
s^{n-i}}),\frac{D^{n}T}{\partial s^{n}}>\notag
\end{equation}
\begin{equation}
+2\sum^{n-1}_{i=0}<\frac{D^{i}}{\partial
s^{i}}(R(T,\frac{DT}{\partial s})\frac{D^{n-1-i}T}{\partial
s^{n-1-i}}),\frac{D^{n}T}{\partial s^{n}}>.\notag
\end{equation}
It is easy to see that \setcounter{equation}{1}
\begin{equation}
\frac{D^{i}}{\partial s^{i}}(k^{2}\frac{D^{n-i}T}{\partial
s^{n-i}})=\sum_{j+k\leq
i}\mathcal{C}_{ijk}<\frac{D^{j+1}T}{\partial
s^{j+1}},\frac{D^{k+1}T}{\partial
s^{k+1}}>\frac{D^{n-j-k}T}{\partial s^{n-j-k}},
\end{equation}
and
\begin{equation}
\frac{D^{i}}{\partial s^{i}}(R(T,\frac{DT}{\partial
s})\frac{D^{n-1-i}T}{\partial s^{n-1-i}})=\sum_{j+k\leq
i}\mathcal{C}_{ijk}R(\frac{D^{j}T}{\partial
s^{j}},\frac{D^{k+1}T}{\partial
s^{k+1}})\frac{D^{n-1-j-k}T}{\partial s^{n-1-j-k}},
\end{equation}
where the coefficients $C_{ijk}$ are constants. To obtain (3), we
have repeatedly used (1).\par Noting that $\mathbf{M}$ satisfies
Condition ($\Lambda$), and then putting above equations together,
we obtain
\begin{equation}
\frac{\partial}{\partial t}(|\frac{D^{n}T}{\partial
s^{n}}|^{2})\leq \frac{\partial^{2}}{\partial
s^{2}}(|\frac{D^{n}T}{\partial
s^{n}}|^{2})-|\frac{D^{n+1}T}{\partial
s^{n+1}}|^{2}+C_{1}|\frac{D^{n}T}{\partial
s^{n}}|^{2}+C_{2}|\frac{DT}{\partial
s}|^{2}|\frac{D^{n}T}{\partial s^{n}}|^{2}\notag
\end{equation}
\begin{equation}
+C_{3}|\frac{D^{2}T}{\partial s^{2}}||\frac{D^{n}T}{\partial
s^{n}}|^{2}+C_{4}\sum_{0\leq i,j,k<n}|\frac{D^{i}T}{\partial
s^{i}}||\frac{D^{j}T}{\partial s^{j}}||\frac{D^{k}T}{\partial
s^{k}}||\frac{D^{n}T}{\partial s^{n}}|,\notag
\end{equation}
where $C_{i}$ are positive constants depending on $n$ and
$\Lambda$. In the last term, the range of indices satisfies in
addition either
$$
i+j+k = n+2
$$
or
$$
i+j+k = n.
$$
\begin{Thm}
Fix $t_{0}\in [0,+\infty)$. Let
$$
M_{t_{0}} =: \max k^{2}(\cdot,t_{0}).
$$
Assume
$$
M_{t_{0}} < +\infty.
$$
Then there exist constants $\tilde{c}_{l}<+\infty$ independent of
$t_{0}$ such that for $t\in (t_{0},t_{0}+\frac{1}{2\Lambda}\log
(1+\frac{\Lambda}{4M_{t_{0}}+\Lambda+1})]$, we have
$$
|\frac{D^{l}T}{\partial s^{l}}|^{2}\leq
\frac{\tilde{c}_{l}M_{t_{0}}}{(t-t_{0})^{l-1}}.
$$
\end{Thm}
\begin{proof}
Without loss of generality, we may assume that $t_{0} = 0$, and
then translate the estimates.\par {\it (1)} For $l = 1$, we have
\begin{equation}
\hspace{-10mm}\frac{\partial}{\partial t}(|\frac{DT}{\partial
s}|^{2})=2<\frac{D}{\partial s}\frac{D}{\partial
t}T+k^{2}\frac{DT}{\partial s}+R(T,\frac{DT}{\partial
s})T,\frac{DT}{\partial s}>\notag
\end{equation}
\begin{equation}
\hspace{9mm}=2<\frac{D^{3}T}{\partial s^{3}},\frac{DT}{\partial
s}>+4k^{4}+2k^{2}R(T,N,T,N)\notag
\end{equation}
\begin{equation}
\hspace{-2mm}\leq \frac{\partial^{2}}{\partial
s^{}}(|\frac{DT}{\partial s}|^{2})-2|\frac{D^{2}T}{\partial
s^{2}}|^{2}+4k^{4}+2\Lambda k^{2}\notag
\end{equation}
It follows from the maximum principle that $M_{t}$ satisfies
\begin{equation}
\log \frac{M_{t}}{\frac{2}{\Lambda}M_{t}+1}-\log
\frac{M_{0}}{\frac{2}{\Lambda}M_{0}+1}\leq 2\Lambda t.\notag
\end{equation}
If
$$
t\leq
\frac{1}{2\Lambda}\log(1+\frac{\Lambda}{4M_{t_{0}}+\Lambda+1}),
$$
then
$$
M_{t} \leq 2M_{0}.
$$
So we may choose $\tilde{c_{1}} = 2$.\par {\it (2)} For $l = 2$,
we have
\begin{equation}
\frac{\partial}{\partial t}(|\frac{D^{2}T}{\partial
s^{2}}|^{2})=2<\frac{D}{\partial s}\frac{D}{\partial
t}\frac{DT}{\partial s}+k^{2}\frac{D^{2}T}{\partial
s^{2}}+R(T,\frac{DT}{\partial s})\frac{DT}{\partial
s},\frac{D^{2}T}{\partial s^{2}}>\notag
\end{equation}
\begin{equation}
\hspace{16mm}=2<\frac{D^{2}}{\partial s^{2}}(\frac{DT}{\partial
t}),\frac{D^{2}T}{\partial s^{2}}>+2<\frac{D}{\partial
s}(k^{2}\frac{DT}{\partial s}),\frac{D^{2}}{\partial s^{2}}>\notag
\end{equation}
\begin{equation}
\hspace{11mm}+2<\frac{D}{\partial s}(R(T,\frac{DT}{\partial
s})T),\frac{D^{2}}{\partial s^{2}}>+2k^{2}|\frac{D^{2}T}{\partial
s^{2}}|^{2}\notag
\end{equation}
\begin{equation}
\hspace{-20mm}+2R(T,\frac{DT}{\partial s},\frac{DT}{\partial
s},\frac{D^{2}T}{\partial s^{2}})\notag
\end{equation}
\begin{equation}
\hspace{5mm}\leq \frac{\partial^{2}}{\partial
s^{2}}(|\frac{D^{2}T}{\partial
s^{2}}|^{2})-|\frac{D^{3}T}{\partial
s^{3}}|^{2}+18|\frac{DT}{\partial s}|^{2}|\frac{D^{2}T}{\partial
s^{2}}|^{2}\notag
\end{equation}
\begin{equation}
\hspace{-20mm}+2\Lambda|\frac{D^{2}T}{\partial
s^{2}}|^{2}+2\Lambda|\frac{DT}{\partial s}|^{4}.\notag
\end{equation}
So
\begin{equation}
\hspace{-72mm}\frac{\partial}{\partial t}(t|\frac{D^{2}T}{\partial
s^{2}}|^{2}+3|\frac{DT}{\partial s}|^{2})\notag
\end{equation}
\begin{equation}
\leq \frac{\partial^{2}}{\partial s^{2}}(t|\frac{D^{2}T}{\partial
s^{2}}|^{2}+3|\frac{DT}{\partial s}|^{2})-t|\frac{D^{3}T}{\partial
s^{3}}|^{2}+[t(18|\frac{DT}{\partial
s}|^{2}+2\Lambda)-5]|\frac{D^{2}T}{\partial s^{2}}|^{2}\notag
\end{equation}
\begin{equation}
\hspace{-55mm}+(2\Lambda t+12)|\frac{DT}{\partial
s}|^{4}+6\Lambda|\frac{DT}{\partial s}|^{2}.\notag
\end{equation}
Since
$$
t\leq \frac{1}{2\Lambda}\log (1+\frac{\Lambda}{4M_{0}+\Lambda
+1})\leq \frac{1}{2(4M_{0}+\Lambda+1)},
$$
we have
\begin{equation}
\frac{\partial}{\partial t}(t|\frac{D^{2}T}{\partial
s^{2}}|^{2}+3|\frac{DT}{\partial s}|^{2})\leq
\frac{\partial^{2}}{\partial s^{2}}(t|\frac{D^{2}T}{\partial
s^{2}}|^{2}+3|\frac{DT}{\partial s}|^{2})+52M_{0}^{2}+12\Lambda
M_{0}.\notag
\end{equation}
Thus it follows that
\begin{equation}
t|\frac{D^{2}T}{\partial s^{2}}|^{2}+3|\frac{DT}{\partial
s}|^{2}\leq 16M_{0},\notag
\end{equation}
and we may conclude on this time interval that
$$
|\frac{D^{2}T}{\partial s^{2}}|^{2}\leq \frac{16M_{0}}{t}.
$$
So we may choose $\tilde{c}_{2} = 16$.\par The induction
hypothesis and repeated usage of the Peter-Paul inequality, i.e.,
$$
ab \leq \epsilon a^{2}+\frac{1}{4\epsilon}b^{2},
$$
allow us to find constants $a_{i}$ and $A$, $B$ on our time
interval such that
\begin{equation}
\frac{\partial}{\partial
t}(\sum_{i=1}^{m}a_{i}t^{i-1}|\frac{D^{i}T}{\partial
s^{i}}|^{2})\leq AM_{0}^{2}+BM_{0}.\notag
\end{equation}
Thus we obtain $\tilde{c}_{m}$ as before.
\end{proof}
Note that these estimates prove the long time existence result.
That is, as long as the curvature remains bounded on time interval
$[0,\alpha)$, one can define a smooth limit for the tangent vector
$T$ at time $\alpha$. Thus, by integrating the tangent vector $T$,
one can obtain a smooth limit curve.
\end{section}
\begin{section}{convergence at inifity}
In this section, we want to prove some convergent results which
will be applied in $\S$ 6. Similar results can be found in
~\cite{G89} where the author dealt with Riemannian surfaces. Those
results may be regarded as a generalization of ours at
2-dimensional case.
\par Let $\omega$ be the maximal existence time of the curve
shortening flow. Throughout this section, we shall assume that
$$
\omega = +\infty,
$$
and
$$
\lim_{t\rightarrow +\infty}L(\gamma_{t})>0.
$$
Then we have
\begin{Lem}
The $L^{2}$-norm of curvature converges to zero as $t\rightarrow
+\infty$.
\end{Lem}
\begin{proof}
It is obvious that the time derivative of length,
$$
-\int k^{2}ds,
$$
must be approaching zero at an $\epsilon$-dense set of
sufficiently large times. So what we need to do is to bound its
time derivative.
\begin{equation}
\frac{\partial}{\partial t}\int k^{2}ds\leq \int
\frac{\partial^{2}}{\partial s^{2}}(|\frac{DT}{\partial
s}|^{2})-2|\frac{D^{2}T}{\partial s^{2}}|^{2}+4k^{4}+2\Lambda
k^{2}ds -\int k^{4}ds\notag
\end{equation}
\begin{equation}
\leq -2\int |\frac{D^{2}T}{\partial s^{2}}|^{2}ds+(3\sup
k^{2}+2\Lambda)\int k^{2}ds.\notag
\end{equation}
But
\begin{equation}
\sup k^{2}\leq (\inf |k|+\int |\frac{D^{2}T}{\partial
s^{2}}|^{2}ds)^{2}\leq \frac{2}{l_{\infty}}\int k^{2}ds+2l_{0}\int
|\frac{D^{2}T}{\partial s^{2}}|^{2}ds\notag
\end{equation}
$\Rightarrow$
\begin{equation}
-2\int |\frac{D^{2}T}{\partial s^{2}}|^{2}ds\leq
-\frac{1}{l_{0}}\sup k^{2}+\frac{2}{l_{\infty}l_{0}}\int
k^{2}ds.\notag
\end{equation}
So
\begin{equation}
\frac{\partial}{\partial t}\int k^{2}ds\leq
(2\Lambda+\frac{2}{l_{\infty}l_{0}})\int k^{2}ds+\sup k^{2}(3\int
k^{2}ds-\frac{1}{l_{0}}).\notag
\end{equation}
Therefore $\int k^{2}ds$ has at most exponential growth when it is
sufficiently small. This implies that it must converge to zero at
$t\rightarrow \infty$.
\end{proof}
\begin{Lem}
$$
\lim_{t\rightarrow \infty}\int |\frac{D^{2}T}{\partial
s^{2}}|^{2}ds=0.
$$
\end{Lem}
\begin{proof}
Suppose not. We need only consider those times when $\int
|\frac{D^{2}T}{\partial s^{2}}|^{2}ds$ is sufficiently greater
than $\int k^{2}ds$. Look at the time derivative of $\int
|\frac{D^{2}T}{\partial s^{2}}|^{2}ds$. Our rules for
differentiating yields:
\begin{equation}
\frac{\partial}{\partial t}\int |\frac{D^{2}T}{\partial
s^{2}}|^{2}ds\leq \int -|\frac{D^{3}T}{\partial
s^{3}}|^{2}+18|\frac{DT}{\partial s}|^{2}|\frac{D^{2}T}{\partial
s^{2}}|^{2}+2\Lambda|\frac{D^{2}T}{\partial
s^{2}}|^{2}+\Lambda^{2}|\frac{DT}{\partial s}|^{2}.\notag
\end{equation}
We will bound the last three terms in this integral by a fraction
of the first.\par Note that
\begin{equation}
\frac{\partial}{\partial s}<\frac{DT}{\partial
s},\frac{DT}{\partial s}>=<\frac{DT}{\partial
s},\frac{D^{3}T}{\partial s^{3}}>+|\frac{D^{2}T}{\partial
s^{2}}|^{2}\notag
\end{equation}
$\Rightarrow$
\begin{equation}
0 = \int <\frac{DT}{\partial s},\frac{D^{3}T}{\partial s^{3}}>ds +
\int |\frac{D^{2}T}{\partial s^{2}}|^{2}ds\notag
\end{equation}
$\Rightarrow$
\begin{equation}
\int |\frac{D^{2}T}{\partial s^{2}}|^{2}d \leq (\int
|\frac{DT}{\partial s}|^{2}ds \cdot \int |\frac{D^{3}T}{\partial
s^{3}}|^{2}ds)^{\frac{1}{2}}.\notag
\end{equation}
If we assume that
$$
\int |\frac{D^{2}T}{\partial s^{2}}|^{2}ds>\alpha \cdot \int
|\frac{DT}{\partial s}|^{2}ds,
$$
then we get
\begin{equation}
\int |\frac{D^{2}T}{\partial s^{2}}|^{2}ds \leq \alpha^{-1} \cdot
\int |\frac{D^{3}T}{\partial s^{3}}|^{2}ds.\notag
\end{equation}
Assume that
$$
\int |\frac{DT}{\partial s}|^{2}ds \leq \epsilon
$$
for some small $\epsilon > 0$. We estimate the second term:
\begin{equation}
\int |\frac{DT}{\partial s}|^{2}|\frac{D^{2}T}{\partial
s^{2}}|^{2}ds \leq \epsilon \cdot \sup |\frac{D^{2}T}{\partial
s^{2}}|^{2}.\notag
\end{equation}
But
\begin{equation}
\hspace{-5mm}\sup |\frac{D^{2}T}{\partial s^{2}}|^{2}\leq \{\inf
|\frac{D^{2}T}{\partial s^{2}}|+\int |\frac{\partial}{\partial s
}(|\frac{D^{2}T}{\partial s^{2}}|)|ds\}^{2}\notag
\end{equation}
\begin{equation}
\hspace{4mm}\leq (\inf |\frac{D^{2}T}{\partial s^{2}}|+\int
|\frac{D^{3}T}{\partial s^{3}}|ds)^{2}\notag
\end{equation}
\begin{equation}
\hspace{17mm}\leq \frac{2}{l_{\infty}}\int |\frac{D^{2}T}{\partial
s^{2}}|^{2}ds+2l_{0}\int |\frac{D^{3}T}{\partial
s^{3}}|^{2}ds\notag
\end{equation}
\begin{equation}
\hspace{2mm}\leq (\frac{2}{\alpha l_{\infty}}+2l_{0})\int
|\frac{D^{3}T}{\partial s^{3}}|^{2}ds.\notag
\end{equation}
Hence
\begin{equation}
\frac{\partial}{\partial t}\int |\frac{D^{2}T}{\partial
s^{2}}|^{2}ds\leq [-1+\epsilon \cdot (\frac{2}{\alpha
l_{\infty}}+2l_{0})+\frac{2\Lambda}{\alpha}+\frac{\Lambda^{2}}{\alpha^{2}}]\int
|\frac{D^{3}T}{\partial s^{3}}|^{2}ds\notag
\end{equation}
\begin{equation}
\hspace{-29mm}\leq -\frac{1}{2}\int |\frac{D^{3}T}{\partial
s^{3}}|^{2}ds\notag
\end{equation}
\begin{equation}
\hspace{-29mm}\leq -\frac{1}{2}\int |\frac{D^{2}T}{\partial
s^{2}}|^{2}ds.\notag
\end{equation}
So, either $\int |\frac{D^{2}T}{\partial s^{2}}|^{2}ds$ decays
exponentially, or it is comparable to $\int |\frac{DT}{\partial
s}|^{2}ds$. In either event, it decreases to zero.
\end{proof}
The following Sobolev inequality is useful, and will be used
repeatedly.
\begin{Lem}
If $\|f\|_{2}\leq C$ and $\|f'\|_{2}\leq C$, then
$$\|f\|_{\infty}\leq(\frac{1}{\sqrt{2\pi}}+\sqrt{2\pi})C,$$
where $\|\cdot\|_{2}$ is the $L_{2}$ norm and $\|\cdot\|_{\infty}$
is the sup norm for functions on $\mathbf{S}^{1}$.
\end{Lem}
Notice $$|\frac{\partial}{\partial s}(|\frac{DT}{\partial
s}|)|^{2}\leq |\frac{D^{2}T}{\partial s^{2}}|^{2}.$$ So, from
Lemma 8, we have
$$\lim_{t\rightarrow +\infty}\int |\frac{\partial}{\partial
s}(|\frac{DT}{\partial s}|)|^{2}ds=0.$$ Then it follows from the
Sobolev inequality that $\sup (|\frac{DT}{\partial s}|)$ decreases
to zero.\par

We deal with the higher derivatives in the same fashion.
Integration and the Holder inequality yield:
\begin{equation}
(\int |\frac{D^{n}T}{\partial s^{n}}|^{2}ds)^{2}\leq \int
|\frac{D^{n-1}T}{\partial s^{n-1}}|^{2}ds\cdot \int
|\frac{D^{n+1}T}{\partial s^{n+1}}|^{2}ds.\notag
\end{equation}\par
We start, knowing that
$$
\int |\frac{D^{n-1}T}{\partial s^{n-1}}|^{2}ds \rightarrow 0,
$$
and that
$$
\sup (|\frac{D^{m}T}{\partial s^{m}}|) \rightarrow 0
$$
for all $m<n-1$. Then, as before, we can show that $\int
|\frac{D^{n}T}{\partial s^{n}}|^{2}ds$ decreases exponentially
when it is much bigger than some linear combination of $\int
|\frac{D^{n-1}T}{\partial s^{n-1}}|^{2}ds$ and terms involving the
lower order derivatives. Therefore
$$
\int |\frac{D^{n}T}{\partial s^{n}}|^{2}ds \rightarrow 0.
$$
And then it easily follows from the inequality
$$
|\frac{\partial}{\partial s}(|\frac{D^{n-1}T}{\partial
s^{n-1}}|)|^{2}\leq |\frac{D^{n}T}{\partial s^{n}}|^{2},
$$
and the Sobolev inequality that
$$
\sup (|\frac{D^{n-1}T}{\partial s^{n-1}}|) \rightarrow 0.
$$
We leave the details as an exercise. Or one can refer to
~\cite{G89}.\par By now, we have proved
\begin{Thm} $\mathbf{M}$ is a compact locally symmetric space satisfying
Condition ($\Lambda$). If the curve shortening flow exists for
infinite time, and
$$
\lim_{t\rightarrow \infty}L(\gamma_{t}) > 0,
$$
then for every $n > 0$,
$$
\lim_{t\rightarrow \infty}\sup (|\frac{D^{n}T}{\partial s^{n}}|) =
0.
$$
Moreover, the limiting curve exists and is a geodesic.
\end{Thm}
\end{section}
\begin{section}{Curve shortening in space forms}

In this section, we shall assume that $\mathbf{M}$ is a
3-dimensional Riemannian manifold with constant sectional
curvature $K$. In this case, the curvature operator $R$ has the
following simple expression, i.e.,
\begin{equation}
R(X_{1},X_{2})X_{3}=K(<X_{1},X_{3}>X_{2}-<X_{2},X_{3}>X_{1})\notag
\end{equation}
for all $X_{i}\in T\mathbf{M}$. In particular,
$$
R(T,N)T=KN
$$
and
$$
R(T,N)N=-KT.
$$\par
For $n=3$, we have the well-known Frenet matrix for a curve
$\gamma$, with the arc-length parameter $s$,
\begin{equation}
\frac{D}{\partial s} \left (
\begin{array}{ll}
    T\\
    N\\
    B\\
\end{array}
\right ) = \left (
    \begin{array}{lll}
    0 & k & 0\\
    -k & 0 & \tau\\
    0 & -\tau & 0\\
    \end{array}
\right ) \left (
\begin{array}{ll}
    T\\
    N\\
    B\\
\end{array}
\right ) .\notag
\end{equation}
With these relations, we can make Lemma 5 more precise.
\begin{Lem}
$$
\nabla_{t}T = \frac{\partial k}{\partial s}N + \tau kB.
$$
\end{Lem}
\begin{proof}
Note
\begin{equation}
\frac{D^{2}T}{\partial s^{2}}=\frac{D}{\partial s}(kN)
=\frac{\partial k}{\partial s}N+k(-kT+\tau B),\notag
\end{equation}
so
\begin{equation}
k^{2}T+\frac{D^{2}T}{\partial s^{2}}=\frac{\partial k}{\partial
s}N+\tau kB.\notag
\end{equation}
Then, by Lemma 5, we get
\begin{equation}
\nabla_{t}T=k^{2}T+\frac{D^{2}T}{\partial s^{2}}=\frac{\partial
k}{\partial s}N+\tau kB.\notag
\end{equation}
\end{proof}
Now, we can compute the evolution of curvature $k$.
\begin{Lem}
$$
\frac{\partial k}{\partial t} = \frac{\partial^{2}k}{\partial
s^{2}} + k^{3} - \tau^{2}k + Kk.
$$
\end{Lem}
\begin{proof}
Note
\begin{equation}
\hspace{-8mm}\nabla_{t}\nabla_{s}T=\nabla_{s}\nabla_{t}T+k^{2}\nabla_{s}T+kR(T,N)T\notag
\end{equation}
\begin{equation}
\hspace{7mm}=\nabla_{s}(\frac{\partial k}{\partial s}N+\tau
kB)+k^{3}N+KkN\notag
\end{equation}
\begin{equation}
\hspace{14mm}=\frac{\partial^{2}k}{\partial s^{2}}N+\frac{\partial
k}{\partial s}(-kT+\tau B)+\frac{\partial}{\partial s}(\tau
k)B\notag
\end{equation}
\begin{equation}
\hspace{3mm}+(\tau k)(-\tau N)+k^{3}N+KkN\notag
\end{equation}
\begin{equation}
\hspace{13mm}=-k\frac{\partial k}{\partial
s}T+(\frac{\partial^{2}k}{\partial
s^{2}}+k^{3}-\tau^{2}k+Kk)N\notag
\end{equation}
\begin{equation}
\hspace{-9mm}+(\tau \frac{\partial k}{\partial
s}+\frac{\partial}{\partial s}(\tau k))B,\notag
\end{equation}
so
\begin{equation}
\hspace{-17mm}\frac{\partial k}{\partial
t}=\frac{\partial}{\partial t}<\nabla_{s}T,N>\notag
\end{equation}
\begin{equation}
\hspace{17mm}=<\nabla_{t}\nabla_{s}T,N>+<\nabla_{s}T,\nabla_{t}N>\notag
\end{equation}
\begin{equation}
\hspace{-1mm}=\frac{\partial^{2}k}{\partial
s^{2}}+k^{3}-\tau^{2}k+Kk.\notag
\end{equation}
\end{proof}
We also need to know the rate at which the unit normal vector to
the curve rotates. This can be got directly from the proof of
Lemma 12.
\begin{Coro}
$$
\nabla_{t}N = -\frac{\partial k}{\partial s}T + (\frac{\partial
\tau}{\partial s} + 2\frac{\tau}{k}\frac{\partial k}{\partial
s})B.
$$
\end{Coro}
\begin{proof}
Note
\begin{equation}
\nabla_{s}T=kN,\notag
\end{equation}
so
\begin{equation}
\nabla_{t}\nabla_{s}T=\nabla_{t}(kN)=\frac{\partial k}{\partial
t}N+k\nabla_{t}N.\notag
\end{equation}
With this equation, and noticing that
$$
<\nabla_{t}N,N> = 0,
$$
we get the following relation from the proof of Lemma 12:
\begin{equation}
k\nabla_{t}N=-k\frac{\partial k}{\partial s}T+(k\frac{\partial
\tau}{\partial s}+2\tau \frac{\partial k}{\partial s})B.\notag
\end{equation}
Multiplying both sides $\frac{1}{k}$, we obtain what we want.
\end{proof}
Now, we can compute the evolution of torsion $\tau$.
\begin{Lem}
$$
\frac{\partial \tau}{\partial t} =
\frac{\partial^{2}\tau}{\partial s^{2}} +
2\frac{1}{k}\frac{\partial k}{\partial s}\frac{\partial
\tau}{\partial s} + 2\frac{\tau}{k}(\frac{\partial^{2}k}{\partial
s^{2}} - \frac{1}{k}(\frac{\partial k}{\partial s})^{2} + k^{3}).
$$
\end{Lem}
\begin{proof}
We have
\begin{equation}
\nabla_{s}N=-kT+\tau B\notag
\end{equation}
$\Rightarrow$ \setcounter{equation}{3}
\begin{equation}
\nabla_{t}\nabla_{s}N = -\frac{\partial k}{\partial t}T -
k\nabla_{t}T + \frac{\partial \tau}{\partial t}B + \tau
\nabla_{t}B.
\end{equation}\par
The left hand side of (4) equals
\begin{equation}
\hspace{-17mm}\nabla_{s}\nabla_{t}N+k^{2}\nabla_{s}N+kR(T,N)N\notag
\end{equation}
\begin{equation}
\hspace{17mm}=\nabla_{s}(-\frac{\partial k}{\partial
s}T+(\frac{\partial \tau}{\partial
s}+2\frac{\tau}{k}\frac{\partial k}{\partial s})B)+k^{2}(-kT+\tau
B)-KkT\notag
\end{equation}
\begin{equation}
\hspace{-6mm}=-\frac{\partial^{2}k}{\partial
s^{2}}T-k\frac{\partial k}{\partial
s}N+(\frac{\partial^{2}\tau}{\partial
s^{2}}+2\frac{\partial}{\partial s}(\frac{\tau}{k}\frac{\partial k
}{\partial s}))B\notag
\end{equation}
\begin{equation}
\hspace{5mm}+(\frac{\partial \tau}{\partial
s}+2\frac{\tau}{k}\frac{\partial k }{\partial s})(-\tau
N)-k^{3}T+\tau k^{2}B-KkT.\notag
\end{equation}
The coefficient of $B$ is
$$
\frac{\partial^{2}\tau}{\partial s^{2}}+2\frac{\partial}{\partial
s}(\frac{\tau}{k}\frac{\partial k }{\partial s})+\tau k^{2}.
$$\par
The right hand side of (4) equals
\begin{equation}
-\frac{\partial k}{\partial t}T-k(\frac{\partial k}{\partial
s}N+\tau kB)+\frac{\partial \tau}{\partial t}B+\tau
\nabla_{t}B.\notag
\end{equation}
Note that
$$
<\nabla_{t}B,B> = 0.
$$
So the coefficient of $B$ is
$$
-\tau k^{2}+\frac{\partial \tau}{\partial t}.
$$\par Hence we get
\begin{equation}
\frac{\partial^{2}\tau}{\partial s^{2}}+2\frac{\partial}{\partial
s}(\frac{\tau}{k}\frac{\partial k }{\partial s})+\tau k^{2}=-\tau
k^{2}+\frac{\partial \tau}{\partial t}\notag
\end{equation}
$\Rightarrow$
\begin{equation} \frac{\partial \tau}{\partial
t}=\frac{\partial^{2}\tau}{\partial
s^{2}}+2\frac{1}{k}\frac{\partial k}{\partial s}\frac{\partial
\tau}{\partial s}+2\frac{\tau}{k}(\frac{\partial^{2}k}{\partial
s^{2}}-\frac{1}{k}(\frac{\partial k}{\partial
s})^{2}+k^{3}).\notag
\end{equation}
\end{proof}
If both $k$ and $\tau$ only depend on time, then their evolutions
become \setcounter{equation}{4}
\begin{equation}
\{
    \begin{array}{ll}
    \frac{dk}{dt}=k^{3}-\tau^{2}k+Kk\\
    \frac{d\tau}{dt}=2\tau k^{2}\\
    \end{array}
\end{equation}

In the following, we study the system in space forms. Since the
flat case $R^3$ was treated by S.J.Altschuler and M.A.Grayson
\cite{AG91}, we consider two remaining cases.

 {\it Case (1): } $\mathbf{M}=\mathbf{H}^{3}$\par In this case,
$K=-1$. The system (5) becomes
\begin{equation}
\{
    \begin{array}{ll}
    \frac{dk}{dt}=k^{3}-k-\tau^{2}k\\
    \frac{d\tau}{dt}=2\tau k^{2}\\
    \end{array}
\end{equation}
We solve (6) as follows.\par For the sake of simplicity, we
introduce the following notations. Let
$$
\{
    \begin{array}{ll}
    u =: k^{2},\\
    v =: \tau^{2}.\\
    \end{array}\notag
$$
Then the system (6) is equivalent to \setcounter{equation}{6}
\begin{equation}
\{
    \begin{array}{ll}
    \frac{du}{dt}=2u^{2}-2u-2uv \cdots (*)\\
    \frac{dv}{dt}=4uv \cdots (**)\\
    \end{array}
\end{equation}\par
Without loss of generality, we may assume the initial condition
$$
v(0)>0.
$$
From $(**)$, we know that $v$ is non-decreasing. So
$$
v(t) > 0
$$
for all $t \geq 0$. Then we can divide $(*)$ by $(**)$ to get
\begin{equation}
\frac{du}{dv}=\frac{u-1}{2v}-\frac{1}{2}.\notag
\end{equation}
Define
$$
w =: u-1,
$$
then \setcounter{equation}{7}
\begin{equation}
\frac{dw}{dv}=\frac{1}{2}(\frac{w}{v}-1).
\end{equation}
Also define
$$
z =: \frac{w}{v},
$$
then
$$
w = zv.
$$
Substituting it into (8), we get
\begin{equation}
\frac{dz}{dv}=\frac{\frac{1}{2}(z-1)-z}{v}=-\frac{z+1}{2v}.\notag
\end{equation}
Integrating from time $0$ to $t$, we have
\begin{equation}
\frac{z(v(t))+1}{z(v(0))+1}=(\frac{v(t)}{v(0)})^{-\frac{1}{2}}.\notag
\end{equation}
Note
$$
z = \frac{u-1}{v}.
$$
So
\begin{equation}
\frac{\frac{u(v(t))-1}{v(t)}+1}{\frac{u(v(0))-1}{v(0)}+1}=(\frac{v(t)}{v(0)})^{-\frac{1}{2}}\notag
\end{equation}
$\Rightarrow$
\begin{equation}
u(v(t))=1+v(t)[(\frac{v(t)}{v(0)})^{-\frac{1}{2}}(\frac{u(v(0))-1}{v(0)}+1)-1]\notag
\end{equation}
\begin{equation}
\hspace{13mm}=1+(v(t)v(0))^{\frac{1}{2}}(\frac{u(v(0))-1}{v(0)}+1)-v(t).\notag
\end{equation}\par
Still, without loss of generality, we may assume another initial
condition
$$
u(v(0))=1.
$$
Then \setcounter{equation}{8}
\begin{equation}
u(v(t))=1+(v(t)v(0))^{\frac{1}{2}}-v(t)
\end{equation}
Substituting (9) into $(**)$, we have
\begin{equation}
\frac{dv}{dt}=4v(1+(vv_{0})^{\frac{1}{2}}-v),\notag
\end{equation}
where
$$
v_{0} =: v(0).
$$\par Let
$$
\tilde{\tau} =: \sqrt{v},
$$
i.e.,
$$
v = \tilde{\tau}^{2},
$$
then \setcounter{equation}{9}
\begin{equation}
\frac{d\tilde{\tau}}{dt}=2\tilde{\tau}(1+\tilde{\tau}_{0}\tilde{\tau}-\tilde{\tau}^{2}),
\end{equation}
where
$$
\tilde{\tau}_{0} =: \tilde{\tau}(0).
$$
Solving (10), we obtain
\begin{equation}
\hspace{-16mm}\tilde{\tau}^{a}(-\tilde{\tau}+\frac{\tilde{\tau}_{0}
+\sqrt{\tilde{\tau}_{0}^{2}+4}}{2})^{b}
(\tilde{\tau}-\frac{\tilde{\tau}_{0}-\sqrt{\tilde{\tau}_{0}^{2}+4}}{2})^{c}\notag
\end{equation}
\begin{equation}
=\tilde{\tau}_{0}^{a}(\frac{-\tilde{\tau}_{0}+\sqrt{\tilde{\tau}_{0}^{2}+4}}{2})^{b}
(\frac{\tilde{\tau}_{0}+\sqrt{\tilde{\tau}_{0}^{2}+4}}{2})^{c}\cdot
\exp(-2t), \cdots \cdots (\diamond)\notag
\end{equation}
where
$$\hspace{-20mm}a=-1,$$
$$b=\frac{1}{2}-\frac{\tilde{\tau}_{0}}{2\sqrt{\tilde{\tau}_{0}^{2}+4}},$$
$$c=\frac{1}{2}+\frac{\tilde{\tau}_{0}}{2\sqrt{\tilde{\tau}_{0}^{2}+4}}.$$
Notice that both $b$ and $c$ are positive. Let $t\rightarrow
+\infty$, then the right hand side of $(\diamond)$ tends to $0$.
So it must be
$$
\tilde{\tau}\rightarrow\frac{\tilde{\tau}_{0}+\sqrt{\tilde{\tau}
_{0}^{2}+4}}{2}.
$$
Together with (9), we see that
$$u \rightarrow 0+.
$$

\par {\it Case (2):}
$\mathbf{M}=\mathbf{S}^{3}$\par In this case, $K=1$. Using the
same method as in {\it (1)}, we know that, as $t\rightarrow
+\infty$,
$$
\sqrt{v}\rightarrow\frac{m+\sqrt{m^{2}+4}}{2},
$$
where
$$
m =: \sqrt{v(0)}\cdot (\frac{2}{v(0)}+1),$$ and
$$u\rightarrow 0+.$$
\begin{Rk}
In both cases, the limiting curves, if they exist, are geodesics.
Moreover, the non-zero torsion reflects the fact that the frames
are twisting along the geodesics.
\end{Rk}
\end{section}

\begin{section}{ramps in the flow}

In this section, we deal with product Riemannian manifolds
$(\mathbf{M}\times \mathbf{S}^{1},g+d\sigma^{2})$. As before,
define
$$
\gamma_{t}(\cdot) =: \gamma(\cdot,t):\mathbf{S}^{1}
\rightarrow\mathbf{M}\times \mathbf{S}^{1}
$$
is an evolving immersed curve along the curve shortening flow. Let
$$
\pi_{\mathbf{S}^{1}}:\mathbf{M}\times\mathbf{S}^{1}\rightarrow
\mathbf{S}^{1}
$$
be projection. It naturally induces a linear mapping
$$
(\pi_{\mathbf{S}^{1}})_{\ast} : T_{\cdot}(\mathbf{M}\times
\mathbf{S}^{1})\rightarrow T_{\pi_{\mathbf{S}^{1}}(\cdot)}S^{1}.
$$
\begin{Def}[definition]
We shall call $\gamma_{t}$ a ramp if there exists a unit tangent
vector field $U$ to $\mathbf{S}^{1}$ such that
$$
<(\pi_{\mathbf{S}^{1}})_{\ast}(T),U>_{S^{1}}>0
$$
along $\gamma_{t}$.
\end{Def}

From this definition, it is easy to deduce the following
\begin{prop}
An immersed curve is a ramp iff the $T\mathbf{S}^{1}$-component of
its tangent vector is non-zero everywhere.
\end{prop}
\begin{Rk}
Ramp is not a new concept. In fact, many authors have studied it
before (see ~\cite{AG92}, ~\cite{EH89}, ~\cite{P03}). As we will
see, ramps have very good properties.
\end{Rk}
{\it Claim:} For a curve shortening flow, if $\gamma_{0}$ is a
ramp, then for all $t>0$, $\gamma_{t}$ is a ramp, too.\par
\begin{proof}
By definition, there exists a unit tangent vector field $U\in
TS^{1}$, such that
$$
u =: <(\pi_{\mathbf{S}^{1}})_{\ast}(T),U>_{S^{1}}>0
$$
at $t=0$. The time derivative of $u$ is \setcounter{equation}{10}
\begin{equation}
\frac{\partial u}{\partial t}=u^{\prime\prime}+k^{2}u.
\end{equation}
Here $\prime$ denotes differential with respect to $s$. If we
define
$$
\mu_{t} =: \min_{S^{1}}u(\cdot,t),
$$
then $\mu_{0}>0$. (11) tells us that $\mu_{t}$ is non-decreasing.
So we obtain the {\it Claim}.
\end{proof}
\begin{prop}
Assume the sectional curvature of $\mathbf{M}\times
\mathbf{S}^{1}$ has an upper bound $\Xi>0$, and $\gamma_{0}$ is a
ramp.\par {\it (1)} Let
$$
\kappa_{t} =: \min_{\mathbf{S}^{1}}k(\cdot,t).
$$
If $\kappa_{t}<0$ for all $t\geq 0$, then
$$
k(\cdot,t)\geq C_{1}\exp(\Xi t)
$$
for all $t\geq 0$, where $C_{1}$ is negative and only depends on
$\gamma_{0}$.\par {\it (2)} Let
$$
\lambda_{t} =: \max_{\mathbf{S}^{1}}k(\cdot,t).
$$
If $\lambda_{t}>0$ for all $t\geq 0$, then
$$
k(\cdot,t)\leq C_{2}\exp(\Xi t)
$$
for all $t\geq 0$, where $C_{2}$ is positive and only depends on
$\gamma_{0}$.
\end{prop}
\begin{proof}
Since $\gamma_{0}$ is a ramp, our {\it Claim} guarantees that
$\gamma_{t}$ is always a ramp. So we may divide $k$ by $u$. The
time derivative of $\frac{k}{u}$ is
\begin{equation}
\frac{\partial}{\partial
t}(\frac{k}{u})=\frac{\partial^{2}}{\partial
s^{2}}(\frac{k}{u})+2\frac{u^{\prime}}{u}\frac{\partial}{\partial
s}(\frac{k}{u})+\frac{k}{u}(k^{2}-|\frac{DN}{\partial
s}|^{2})+\frac{k}{u}R(T,N,T,N).
\end{equation}
Notice the third term of the right hand side of (12) is
non-positive, and the sectional curvature $R(T,N,T,N)$ is bounded
from above by $\Xi$.\par {\it (1)} If we define
$$
\Phi_{t} =: \min_{\mathbf{S}^{1}}\frac{k}{u}(\cdot,t),
$$
then, by the assumption,
$$
\Phi_{t} < 0
$$
for all $t\geq 0$. Formula (12) tells us that $\Phi_{t}$ satisfies
\begin{equation}
\frac{\partial}{\partial t}\Phi_{t}\geq
\Phi_{t}(k^{2}-|\frac{DN}{\partial
s}|^{2})+\Phi_{t}R(T,N,T,N)\notag
\end{equation}
$\Rightarrow$
\begin{equation}
\Phi_{t}^{-1}\frac{\partial}{\partial t}\Phi_{t} \leq \Xi\notag
\end{equation}
$\Rightarrow$
\begin{equation} \Phi_{t} \geq \Phi_{0}\exp(\Xi
t).\notag
\end{equation}
Note $\Phi_{0}<0$ and $u\leq 1$. Then we can obtain {\it (1)}
easily.\par {\it (2)} If we define
$$
\Psi_{t} =: \max_{\mathbf{S}^{1}}\frac{k}{u}(\cdot,t),
$$
then, by the assumption,
$$
\Psi_{t} > 0
$$
for all $t\geq 0$. Formula (12) tells us that $\Psi_{t}$ satisfies
\begin{equation}
\frac{\partial}{\partial t}\Psi_{t}\leq
\Psi_{t}(k^{2}-|\frac{DN}{\partial
s}|^{2})+\Psi_{t}R(T,N,T,N)\notag
\end{equation}
$\Rightarrow$
\begin{equation}
\Psi_{t}^{-1}\frac{\partial}{\partial t}\Psi_{t} \leq \Xi\notag
\end{equation}
$\Rightarrow$
\begin{equation}
\Psi_{t} \leq \Psi_{0}\exp(\Xi t).\notag
\end{equation}
Note $\Psi_{0}>0$ and $u\leq 1$. Then we can obtain {\it (2)}
easily.\par By now, Proposition 18 is proved.
\end{proof}
The following theorem is a direct consequence of Theorem 10 and
Proposition 18.
\begin{Thm}
$\mathbf{M}\times \mathbf{S}^{1}$ is a compact locally symmetric
space. If $\gamma_{0}$ is a ramp, then the curve shortening flow
will converge to a geodesic in the $C^{\infty}$ norm.
\end{Thm}
\begin{proof}
For $\gamma_{0}$ is a ramp, Proposition 18 guarantees that the
curve shortening flow will not blow-up in finite time. This means
that the flow will exists for infinite time. Moreover, from the
proof of {\it Claim}, we know that $\mu_{t}$ is non-decreasing.
This will guarantee that
$$
\lim_{t\rightarrow \infty}L(\gamma_{t}) > 0.
$$
Then Theorem 10 tells us that the limiting curve exists and is a
geodesic.
\end{proof}
\begin{Rk}
As we know, closed geodesics theory is a fundamental part of
Riemannian geometry. There are a lot of nice works in this theory
(~\cite{K}). In 1929, L.Lusternik and L.Schnirelmann
(~\cite{LS29}) outlined that any Riemannian 2-sphere has at least
three simple closed geodesics. Unfortunately, there was a
shortcoming in their proof. Later, M.Grayson used the curve
shortening flow to prove the three geodesics theorem beautifully
(~\cite{G89}). But, for higher dimensional case, there are only
few results. The importance of Theorem 19 is that it points out a
possible way to find closed geodesics on a compact locally
symmetric space $\mathbf{M}\times \mathbf{S}^{1}$, i.e., finding
closed ramps representing non-trivial homology classes in the path
space $(\Sigma,\Sigma_{0})$ of closed curves relative to the point
curve, and then evolving them along the curve shortening flow.
\end{Rk}
\end{section}

\begin{section}{shortening curves in evolving metric}
In this section, we want to study how to evolve the metric to make
a specific curve shortening flow exists as long as the metric is
non-singular. One motivation for this problem is G.Perelman's work
(~\cite{P03}) where he considers the curve shortening flow during
the Ricci flow.\par First, we consider a simple case.

\par Let $\mathbf{M}$ be an oriented
differentiable manifold. Let $g_{0}$ is a metric on $M$ at $t=0$.
For an initial curve $\gamma_{0}$, we shall evolve it according to
the curve shortening flow. Meanwhile, we shall change the metric
as time goes on according to the conformal flow. More precisely,
let
$$
g_{t} = \exp(f)g_{0}
$$
be the time-dependent metric on $M$ at time $t \geq 0$, where
$$
f:M\times [0,+\infty)\rightarrow \mathbf{R}
$$
is a $C^{\infty}$ function, satisfying
$$
f(\cdot,0) \equiv 0.
$$
The following Theorem 21 tells us that if we evolve $f$ properly,
the flow will never blow-up if the metric keeps non-singular.
\begin{Thm}
Define
$$
\Lambda_{t} =: \{p\in \mathbf{S}^{1}|k^{2}(p,t)=M_{t}\}.
$$
If $f$ satisfies
$$
\frac{\partial f}{\partial t}\mid_{\Lambda_{t}} = 2k^{2} +
2R(T,N,T,N),
$$
then the curve shortening flow will exist as long as $g_{t}$ keeps
non-singular.
\end{Thm}
\begin{proof}
The different from before is that the metric on $\mathbf{M}$
depends on time. So we need to consider its time derivative. For
convenience, we will use $<\cdot,\cdot>_{t}$ representing $g_{t}$.
Let
$$
v_{t} =: (<\frac{\partial \gamma}{\partial u},\frac{\partial
\gamma}{\partial u}>_{t})^{\frac{1}{2}}.
$$
Then the time derivative of $v_{t}$ is
\begin{equation}
\frac{\partial v_{t}}{\partial
t}=-k^{2}v_{t}+\frac{1}{2}\frac{\partial f}{\partial
t}v_{t}.\notag
\end{equation}
The second term comes from the change of metric (compare this
equation with Lemma 3).\par Next we want to show that covariant
differentiation with respect to $s$ and $t$ are related by the
equation:\par
\begin{equation}
\nabla_{t}\nabla_{s}=\nabla_{s}\nabla_{t}+(k^{2}-\frac{1}{2}\frac{\partial
f}{\partial t})\nabla_{s}+R(T,\frac{DT}{\partial s}).\notag
\end{equation}
In fact,
\begin{equation}
\hspace{-5mm}\nabla_{t}\nabla_{s}=\frac{\partial}{\partial
t}(\frac{1}{v_{t}})\nabla_{u}+\frac{1}{v_{t}}(\nabla_{u}\nabla_{t}
+R(\frac{\partial}{\partial u},\frac{\partial}{\partial t}))\notag
\end{equation}
\begin{equation}
\hspace{5mm}=(k^{2}-\frac{1}{2}\frac{\partial f}{\partial
t})\nabla_{s}+\nabla_{s}\nabla_{t}+R(T,\frac{DT}{\partial
s})\notag
\end{equation}
\begin{equation}
\hspace{5mm}=\nabla_{s}\nabla_{t}+(k^{2}-\frac{1}{2}\frac{\partial
f}{\partial t})\nabla_{s}+R(T,\frac{DT}{\partial s}).\notag
\end{equation}\par
Last, the covariant differentiation of $T$ with respect to time
$t$ is
\begin{equation}
\nabla_{t}T=\frac{\partial}{\partial
t}(\frac{1}{v_{t}})\frac{\partial \gamma}{\partial
u}+\frac{1}{v_{t}}\nabla_{u}\frac{\partial \gamma}{\partial
t}=(k^{2}-\frac{1}{2}\frac{\partial f}{\partial
t})T+\frac{D^{2}T}{\partial s^{2}}.\notag
\end{equation}\par
Now, with these preparations, we can compute the evolution of
$k^{2}$. The desired result is
\begin{equation}
\frac{\partial}{\partial
t}(k^{2})=(k^{2})^{\prime\prime}-2(|\frac{D^{2}T}{\partial
s^{2}}|^{2}-k^{4})+k^{2}[2k^{2}+2R(T,N,T,N)-\frac{\partial
f}{\partial t}].\notag
\end{equation}
Note that the second term on right hand side is non-positive. So
if $f$ satisfy the evolution equation
\begin{equation}
\frac{\partial f}{\partial t} \mid_{\Lambda_{t}} = 2k^{2} +
2R(T,N,T,N),\notag
\end{equation}
then $M_{t}$ is non-increasing. Therefore, it is bounded by
$M_{0}$ for all $t>0$. This implies that, at time $t$, the curve
can go on flowing if $g_{t}$ is non-singular.
\end{proof}
\begin{Ex}
Consider on $\mathbf{R}^{2}$,
$$
g_{0}=dx^{2}+dy^{2}.
$$
Let
$$
\gamma_{0}(\theta) =: (\cos\theta,\sin\theta), \theta \in
\mathbf{S}^{1},
$$
i.e., the unit circle. For this initial curve, we have
$$
k^{2}(\cdot,0) \equiv 1.
$$
\par If we
change the metric according to the conformal equation
$$
g_{t}=\exp(2t)g_{0},
$$
i.e., let
$$
f \equiv 2t,
$$
then the evolving curve satisfies
$$
\gamma_{t}=\exp(-t)\gamma_{0}.
$$
It is easy to see that
$$
k^{2}(\cdot,t) \equiv 1
$$
for all $t > 0$.
\par
Observe that the circle collapses to a point in infinite time
rather than finite time. Similar phenomenon also appears in
~\cite{AL86} where the authors modified the usual curve shortening
flow and introduced new time parameter. We will not go further in
that direction. Readers who are interested in this topic can read
their excellent paper.
\end{Ex}
Finally, we will investigate a little more complecated case.\par
$(\mathbf{M}\times\mathbf{S}^{1},g+\exp(f)d\sigma^{2})$ is a
warped Riemannian manifold (for a warped Riemannian manifold, see,
for example, ~\cite{M03}). In this case, we will fix $g$, and only
change $f$ as time goes on. We want to decide how to change $f$ to
make a specific curve shortening flow exist as long as the metric
keeps non-singular.\par As before, let
$$
\pi_{\mathbf{S}^{1}} : \mathbf{M}\times \mathbf{S}^{1}\rightarrow
\mathbf{S}^{1}
$$
be projection, with induced mapping
$$
(\pi_{\mathbf{S}^{1}})_{\ast} : T_{\cdot}\mathbf{M}\times
\mathbf{S}^{1}\rightarrow
T_{\pi_{\mathbf{S}^{1}}(\cdot)}\mathbf{S}^{1}.
$$
Still, for convenience, we will use $<\cdot,\cdot>_{t}$ and
$<\cdot,\cdot>_{t}^{\mathbf{S}^{1}}$ representing metrics on
$\mathbf{M}\times \mathbf{S}^{1}$ and $\mathbf{S}^{1}$
respectively. Let
$$
v_{t} =: (<\frac{\partial \gamma}{\partial u},\frac{\partial
\gamma}{\partial u}>_{t})^{\frac{1}{2}}.
$$
Then the time derivative of $v_{t}$ is
\begin{equation}
\frac{\partial v_{t}}{\partial
t}=-k^{2}v_{t}+\frac{1}{2}\frac{\partial f}{\partial
t}(|(\pi_{\mathbf{S}^{1}})_{\ast}(T)|_{t}^{\mathbf{S}^{1}})^{2}v_{t}.\notag
\end{equation}
A straightforward calculation shows
\begin{equation}
\nabla_{t}\nabla_{s}=\nabla_{s}\nabla_{t}+(k^{2}-\frac{1}{2}\frac{\partial
f}{\partial
t}(|(\pi_{\mathbf{S}^{1}})_{\ast}(T)|_{t}^{\mathbf{S}^{1}})^{2})\nabla_{s}
+R(T,\frac{DT}{\partial s}),\notag
\end{equation}
and
\begin{equation}
\nabla_{t}T=(k^{2}-\frac{1}{2}\frac{\partial f}{\partial
t}(|(\pi_{\mathbf{S}^{1}})_{\ast}(T)|_{t}^{\mathbf{S}^{1}})^{2})T+\frac{D^{2}T}{\partial
s^{2}}.\notag
\end{equation}
Now we compute the evolution of $k^{2}$:
\begin{equation}
\frac{\partial}{\partial
t}(k^{2})=(k^{2})^{\prime\prime}-2(|\frac{D^{2}T}{\partial
s^{2}}|^{2}-k^{4})+k^{2}\{2k^{2}+2R(T,N,T,N)\notag
\end{equation}
\begin{equation}
\hspace{5mm}-\frac{\partial f}{\partial
t}[2(|(\pi_{\mathbf{S}^{1}})_{\ast}(T)|_{t}^{\mathbf{S}^{1}})^{2}
-(|(\pi_{\mathbf{S}^{1}})_{\ast}(N)|_{t}^{\mathbf{S}^{1}})^{2}]\}.\notag
\end{equation}
If
$$
2(|(\pi_{\mathbf{S}^{1}})_{\ast}(T)|_{t}^{\mathbf{S}^{1}})^{2}
-(|(\pi_{\mathbf{S}^{1}})_{\ast}(N)|_{t}^{\mathbf{S}^{1}})^{2}
\neq 0
$$
on $\Lambda_{t}$ at time $t$, then we can require $f$ satisfy
\begin{equation}
\frac{\partial f}{\partial
t}\mid_{\Lambda_{t}}=\frac{2k^{2}+2R(T,N,T,N)}{2(|(\pi_{\mathbf{S}^{1}})
_{\ast}(T)|_{t}^{\mathbf{S}^{1}})^{2}
-(|(\pi_{\mathbf{S}^{1}})_{\ast}(N)|_{t}^{\mathbf{S}^{1}})^{2}}.\notag
\end{equation}
So $M_{t}$ is non-increasing. By now, we have eventually proved
\begin{Thm}
Assume
$$
2(|(\pi_{\mathbf{S}^{1}})_{\ast}(T)|_{t}^{\mathbf{S}^{1}})^{2}
-(|(\pi_{\mathbf{S}^{1}})_{\ast}(N)|_{t}^{\mathbf{S}^{1}})^{2}
\neq 0
$$
on $\Lambda_{t}$ for all $t\geq 0,$ and $f$ satisfies
$$
\frac{\partial f}{\partial t}
\mid_{\Lambda_{t}}=\frac{2k^{2}+2R(T,N,T,N)}{2(|(\pi_{\mathbf{S}^{1}})
_{\ast}(T)|_{t}^{\mathbf{S}^{1}})^{2}
-(|(\pi_{\mathbf{S}^{1}})_{\ast}(N)|_{t}^{\mathbf{S}^{1}})^{2}}.
$$ Then the
curve shortening flow will exist as long as $g_{t}$ keeps
non-singular.
\end{Thm}
\end{section}

\end{document}